\documentclass[english, 12pt]{article}
\usepackage{amssymb, amsmath, amsfonts}
\usepackage{url}
\usepackage{hyperref}

\newtheorem{definition}{Definition}[section]
\newtheorem{proposition}[definition]{Proposition}

\newtheorem{theorem}[definition]{Theorem}

\begin{document}

\title{Bounds for the maximum modulus of polynomial roots  with nearly optimal worst-case overestimation }
\author{Prashant Batra\\
	{\small Hamburg University of Technology,  D-21071 Hamburg.}\\
	{\small e-mail: batra@tuhh.de, Phone: ++49(40)42878-3478.}\\
	\url{https://orcid.org/0000-0002-4079-3792} }

\date{\today}
\maketitle

\begin{abstract}
	Many  upper bounds for the moduli of polynomial roots have been proposed but  reportedly  assessed on selected examples or restricted classes only. Regarding  quality measured in terms of worst-case relative overestimation of the maximum root-modulus we establish a simple, nearly optimal result.
\end{abstract}

\section{Introduction }
The general solution of polynomial equations via algebraic expressions is impossible  according to the Ruffini-Abel theorem  \cite{BurnsidePanton,McNameeBookvolumeTWO}. Over time a multitude of numerical methods to approximate solutions \cite{SendovHandbook,McNameeBookPartOne,McNameeBookvolumeTWO} or to estimate the root moduli \cite{RahmanSchmeisser_Book, Marden} has evolved. A special mention is warranted for Kalantari's  infinite family of modulus bounds \cite{Kalantari2005Family} which has been shown by Jin \cite{Jin2006Sharpness} to converge to the extremal root-modulus, see also  \cite{HsuCheng2014GRaeffeKalantari} for reference.  
We want to consider in this note upper bounds for the largest modulus of roots of
$p(z)=a_nz^n+a_{n-1}z^{n-1}+ \ldots +a_1z+a_0$. The bounds should be explicit algebraic expressions of the absolute values $|a_i|.$ Thus, we will be dealing with \emph{a priori} bounds  ideally  of low computational effort.   

New inventive methods  to approximate the largest root-modulus via algebraic expressions are still being developed, e.g. \cite{Kittaneh2020Filomat, Melman2013LAMA,Melman2024Cauchy}, but the results are seldom assessed in a fully choice-free fashion. The usual assessment of any bound is carried out through individual evaluation for carefully selected polynomials or via batch-testing sets of polynomials with a given coefficient distribution. We want to consider in this note a generally applicable quality measure for modulus bounds established by van der Sluis \cite{vdSl70} in 1970. This measure avoids the necessity of any choices of examples or coefficient distributions. We  recall van der Sluis' fundamental results in the subsequent  Section \ref{SEC-KnownResults}. We derive a new bound based on the parametrized Cassini ovals 
in Section \ref{SEC-NewBounds}. In Section \ref{SEC_Quality} we show that the parameter choice in our bound leads to a nearly optimal result.
\subsection{Known results and thresholds}\label{SEC-KnownResults}

We denote the maximum modulus of the roots of non-constant $ p \in \mathbb{C}[z]$ by $\mu(p) := \max \{ |\lambda|:\lambda \in \mathbb{C}, p(\lambda)=0\}.$ The   \emph{relative} overestimation of the maximum root- modulus by a given root bound $B(p)$ is defined (\emph{cf.} \cite{vdSl70}) as $B(p)/\mu(p)$ for monic $p \not \equiv z^n$.  
An upper bound for the largest root-modulus  employing only the absolute values of the coefficients is called an \emph{absolute root bound} here (compare \cite{vdSl70}). 
 It was shown by van der Sluis (\cite{vdSl70}, Th.3.8)  that the \emph{Cauchy bound} $\rho(p)$   of a \emph{monic} polynomial
 \begin{align}\label{DefMonicPolynomial}
 p(z)=z^n + \sum_{i=0}^{n-1} a_iz^i,
 \end{align} defined as the largest non-negative root of $z^n - \sum_{i=0}^{n-1} |a_i|z^i $ (see also \cite{RahmanSchmeisser_Book}, Def. 8.1.2.),  has smallest worst-case relative overestimation of $\mu(p)$ among all  absolute root bounds. Taking the supremum of the quotients $\rho(p)/\mu(p)$ for non-trivial polynomials of fixed degree $n \geq 1$ it was shown (cf. \cite{vdSl70},  Th.3.8(e)) that
 \begin{equation}\label{OPTIMUM}
 \limsup_{\stackrel{ p \in \mathbb{C}[z],p(z) \not\equiv z^n}{ p \mbox{ \tiny monic}, \deg(p)=n \geq 1.   } }  \frac{\rho(p)}{\mu(p)}  = 1/(\sqrt[n]{2}-1)\approx 1.442 n \sim n/\log(2).
 \end{equation}
   Thus, no bound using only the coefficient moduli can have a worst-case relative overestimation smaller than $(\sqrt[n]{2}-1)^{-1}$. This remark applies, for example, to all bounds in \cite{RahmanSchmeisser_Book}, Th.8.1.7 and Cor. 8.1.8 as well as  42 of the 45 bounds in   (\cite{McNameeBookPartOne}, Chap.1, pp.28\emph{ff.}) with the exemption of A5, A6 and B5. While  van der Sluis' result also applies, e.g., to the bound of Th.3.3 in \cite{MelmanProcAMS2010} it is not applicable to the explicitly computed, more intricate inclusion sets in \cite{MelmanProcAMS2010} (defined by Cassini ovals with boundary curves of 8th degree) from which said theorem is derived. Unfortunately, the worst-case relative overestimation of available bounds often is considerably larger than $1.442 n$ as can be verified for chosen bounds following \cite{vdSl70} or the arguments leading to Theorem \ref{QualityNewBound} below.
   
    There is no explicit \emph{a priori} coefficient expression of the Cauchy bound $\rho(p).$ Hence, we face the problem to produce an approximation with low worst-case relative overestimation via \emph{a priori} calculated expressions.   For a monic $p$ given by \eqref{DefMonicPolynomial} the Fujiwara bound $F(p)$ (see \cite{Fujiwara}, or, e.g., \cite{RahmanSchmeisser_Book}, Theorem 8.1.7 (ii), or, \cite{Marden}, Ch.30, ex.5), is defined by \\ $F(p):=2 \max \{ |a_{n-1}|,\sqrt[2]{|a_{n-2}|},\sqrt[3]{|a_{n-3}|}, \ldots,\sqrt[n-1]{|a_{1}|}; \sqrt[n]{|a_{0}|/2}  \}.$\\
    \noindent  This bound  has worst-case relative overestimation of $2n$ (which estimate is attainable for every $n \geq 3$) \emph{cf.} \cite{vdSl70}. 
  
  
   Modifying the related Lagrange bound,  an  improvement to a bound with worst-case relative overestimation bounded by $1.58n$ was  obtained in \cite{BatraMignotteStefanescu}.
 Instead of small  improvements of  existing upper bounds (via modifications and  case distinctions) 
 we consider here a different approach.   Our technique eventually leads to a bound which is within two percent of the theoretical optimum \eqref{OPTIMUM}.\\
 
 \section{New inclusion circle via Cassini ovals}\label{SEC-NewBounds}
 For a monic polynomial  $p$ of the form \eqref{DefMonicPolynomial} let us consider the classical Frobenius companion matrix $C_F(p)$   with non-trivial last column. Thus, the matrix $C_F=C_F(p)$ has  subdiagonal equal to 1, and $(-a_0,-a_1, \ldots, -a_{n-2},-a_{n-1})^T$ makes up the  last column. A  similarity transform    with the diagonal matrix $S=S(t)=diag(1,t, \ldots,t^{n-1}),$ $ t>0,$ yields
 $S  C_F S^{-1} =:C(t) $ such that 

\begin{align*}  C(t) =
\begin{pmatrix}   0&\ldots & &\ldots & & 0 & -a_{0}/t^{n-1} \\
t&0 & 0&\ldots & & 0 & -a_{1}/t^{n-2} \\
0 & t & 0& \ldots & & 0 & -a_2/t^{n-3}\\
\vdots & & \ddots& & &  & \\
0 & 	&  &            & t &0 & -a_{n-2}/t\\
0 & 	&  &           & 0 &     t & -a_{n-1}
\end{pmatrix}=: (c_{ij}(t))_{i,j=1}^n.
\end{align*}
It is very well-known that the roots of $p$ are the eigenvalues of $C(t)$, see, e.g., \cite{Marden}.
To locate the eigenvalues of $C(t)$ we employ the reduced row sums $ r_i(t) := \sum_{k \neq i} |c_{ik}(t)| .$ 
It is well-known (see, e.g., \cite{VargaBookGershgorin})  that the eigenvalues of $C(t)$ are contained in the union of  Cassini ovals (the Ostrowski-Brauer sets)
$$ O_{i,j}(t):=\{z \in \mathbb{C}:  |z-c_{ii}(t)|\cdot |z-c_{jj}(t)| \leq r_i(t) \cdot r_j(t) \};i \neq j, i,j=1, \ldots,n .$$ The following bound will be shown in the next section to  have a worst-case relative overestimation close to the optimum stated in \eqref{OPTIMUM}.

%

\begin{proposition}
	Given a complex polynomial $p$ of degree $n \geq 3$ with Taylor expansion $p(z)= \sum_{i=0}^n a_i z^i $, normalized to be monic $(a_n=1)$ and with largest root-modulus $\mu(p).$
	With 
	\begin{equation}\label{FirstParameter}
	{\tau} := \max \{ \sqrt[3]{|a_{n-3}| /2.15}  ;\sqrt[4]{|a_{n-4}|/2}, \sqrt[5]{|a_{n-5}|}, \ldots, \sqrt[n]{|a_{0}|}\}
	\end{equation}  we have that $\mu(p) \leq \Gamma(p) ,$ where $\Gamma(p) $ is defined as
	$$ \max \{ \sqrt{3.15}\sqrt{ {{\tau}}^2+\max \{|a_{n-2}|, 2\tau^2\}}; \frac{|a_{n-1}| + \sqrt{|a_{n-1}|^2 + 4(\tau^2 +  \max \{|a_{n-2}|; 2.15 \tau^2 \}) }}{2}   \}.$$
\end{proposition}

{\bf  Proof:} Let us assume first that $\tau >0,$ and put $t = \tau.$
To estimate the maximum modulus of the roots of $p$ it suffices to bound the points of largest modulus in the Cassini ovals $ O_{i,j}(t)$ for $t=\tau>0.$ In the following, we will repeatedly use the estimate 
\begin{eqnarray}\label{Scaledradius}
(1+|a_{k-1}|/t^{n+1-k}) = (1+|a_{k-1}|/\tau^{n+1-k})  \leq  \begin{cases} 2  & \mbox{for } k \leq n -4; \\
													  (1+2) & \mbox{for } k = n-3 ; \\
													   (1+2.15) &\mbox{for } k = n-2 .\end{cases}
\end{eqnarray}

 \indent 1.) If  $1 \leq j < i \leq n-2$ the Cassini oval $O_{i,j}(t)$ is actually a circular disk around the origin. The radius  $ \sqrt{r_i(\tau)r_j(\tau)}$ equals $ \tau \cdot  \sqrt{(1+|a_{i-1}|/\tau^{n+1-i})(1+|a_{j-1}|/\tau^{n+1-j})}. $
We estimate  the parentheses via \eqref{Scaledradius}, and  find that in the case $i \leq n-2$  all the ovals lie inside
$$ |z| \leq \sqrt{(1+2) \cdot(1+2.15)} \cdot {\tau}.$$

 \indent 2.) Further, if $i=n-1$, $1 \leq j \leq n-2,$ we  write
 $$r_{n-1}(t) \cdot r_j(t)=(t+|a_{n-2}|/t)(t+|a_{j-1}|/t^{n-j})=(t^2+|a_{n-2}|)(1+\frac{|a_{j-1}|}{t^{n+1-j}}).$$ As in the preceding case, the term  $(1+|a_{j-1}|/{\tau}^{n+1-j})$ is bounded for $j \leq n-4$ by $ 2,$   and by $(1+\max \{2;2.15\})$ if $j\geq n-3$. Thus,  the ovals $O_{i,j}(\tau)$ in this case lie in
$$ |z| \leq  \sqrt{ (1+2.15) \cdot ({{\tau}}^2+|a_{n-2}|)}= \sqrt{3.15}\sqrt{  {{\tau}}^2+|a_{n-2}|}.$$  

For the remaining two cases let us first  note that if  $i=n>j$   we have Cassini ovals of the 
form
$$|z+a_{n-1}||z| \leq t(t+a_{j-1}/t^{n-j})=t^2(1+|a_{j-1}|/t^{n+1-j})=:g_{n,j}(t).$$
The point farthest from the origin  lies at a distance\\ $$  \frac{1}{2}(|a_{n-1}| + \sqrt{|a_{n-1}|^2 + 4g_{n,j}(t)  }),$$ 
(see \cite{Melman2013LAMA}, p.186, Sec.2.1). 

3.) For $n-2 \geq j \geq n-3$ and $t=\tau$ we estimate the product $g_{n,j}(t)=g_{n,j}(\tau)$ by  $\tau^2(1+|a_{j-1}|/\tau^{n+1-j})\leq \tau^2(1+\max \{2;2.15\})= 3.15 \tau^2 =: \gamma_{n,j}(\tau).$ 
For $ 1 \leq j \leq n-4,$ we obtain from \eqref{Scaledradius}, the inequality  $g_{n,j}(\tau)=[\tau^2(1+|a_{j-1}|/\tau^{n+1-j})] \leq 2 \tau^2=:\gamma_{n,j}(\tau).$
The Cassini ovals $O_{n,j}$ for $ j \leq n-2$ are thus contained in the circle 
$$ |z| \leq \frac{1}{2}(|a_{n-1}| + \sqrt{|a_{n-1}|^2 + 4 \cdot 3.15 \tau^2  }).$$


4.) The last  oval to  consider stems from the last two rows with $i=n=j+1$, and is given as 
$$O_{n,n-1}( \tau) = \{z \in \mathbb{C}: |z+a_{n-1}||z| \leq \tau^2+|a_{n-2}| \}.$$ 
   This oval is contained in the circle 
$$ |z| \leq \frac{1}{2}(|a_{n-1}| + \sqrt{|a_{n-1}|^2 + 4( \tau^2 + |a_{n-2}|) }).$$

The union of all our circular inclusion regions yield a new root-modulus bound $\Gamma(p)$ for a monic polynomial $p$ whenever $\tau  >0$, namely,
\begin{align*}\Gamma(p) &:=  \max \{ \sqrt{3\cdot 3.15} \cdot \tau; \sqrt{3.15}\sqrt{ {{\tau}}^2+|a_{n-2}|};&\\ & \qquad \qquad \frac{|a_{n-1}| + \sqrt{|a_{n-1}|^2 + 12.6 \tau^2  }}{2}; \frac{|a_{n-1}| + \sqrt{|a_{n-1}|^2 + 4\tau^2 + 4|a_{n-2}| }}{2} \}. &
\end{align*}

 If $\tau$ is equal to zero, then we are essentially dealing with a quadratic polynomial multiplied into $x^{n-2}.$ Trivially, the value $\Gamma(p)$ is a valid upper bound for the root-modulus even if $\tau =0.$
 $\Box$

\section{Relative overestimations of the new root bound}\label{SEC_Quality}
\subsection{Overestimation of the maximum root-modulus}
Our bound functional $\Gamma(\cdot)$ has the following quality.

\begin{theorem}\label{QualityNewBound} For any  $p\in \mathbb{C}[z]$ of degree  $n \geq 3$  s.t. $p(z)=z^n+\sum_{i=0}^{n-1} a_iz^i~\not \equiv~ z^n, $ 
the maximum  relative overestimation of  $\mu(p):=\max \{|\lambda|: p(\lambda)=0\}$ by  $\Gamma(p)$ does not exceed $ 1.4655 n.$ For $n \geq 11$ the overestimation exceeds by less than $5$ per cent the lower threshold \eqref{OPTIMUM} for any absolute root bound.  For $n \geq 85$ the worst-case overestimation realized by our bound  $\Gamma(p)$ exceeds the  lower threshold for any overestimation by at most $2$ per cent. 
\end{theorem}

{\bf Proof:}   Using Vi{\`e}te's representation of the coefficients $a_{n-k}$ as the sum of all possible  products of $k$ different roots we obtain the trivial estimates
$$|a_{n-k}| \leq  \mu(p)^k \binom{n}{k} \leq \mu^k n^k/k! \, ,$$
where we write $\mu$ instead of $\mu(p)$ to save clutter.  
The preceding inequalities imply $\sqrt[k]{|a_{n-k}|}/(n \cdot \mu) \leq \sqrt[k]{1/k!} \, .$
For  $k \geq 5,$ let $c_k$ be $c_k:=c_5:=1/\sqrt[5]{120} \sim  0.4518$, a value  larger or equal to $\sqrt[k]{1/k!} , $ and let $c_4 := \sqrt[4]{1/48}\sim 0.3799, c_3 := \sqrt[3]{1/12.9}\sim 0.4264,c_2:=1/\sqrt{2}\sim 0.7071,c_1:=1.$
With   $ \phi:=\tau/(n \mu)$ (where $\tau$ is defined in \eqref{FirstParameter}) we have
$$ \phi \leq \max \{c_3;c_4,c_5, \ldots,c_n\} = c_3 \sim 0.4264,  \mbox{and}$$
$$ \frac{\Gamma(p)}{n \mu(p)} \leq \max \{\sqrt{9.45} \phi; \sqrt{3.15}\sqrt{ {{\phi}}^2+ 1/2}; \frac{1 + \sqrt{1 + 12.6 \phi^2  } }{2} ;\frac{1 + \sqrt{1 + 4 \phi^2 + 2 }}{2} \}.  $$
  The right-hand side of the preceding inequality  is determined  as (approximately) $\max \{1.3107; 1.4655,1.4070; 1.4653 \}=1.4655.$ The quality claim now follows from  \eqref{OPTIMUM}. $\Box$\\


\subsection{Overestimation of the Cauchy bound}
To assess an absolute root bound fully, van der Sluis  \cite{vdSl70} made additional comparison to the Cauchy bound.
Our new bound $\Gamma(p)$  is homogeneous like the Cauchy bound (i.e., it scales with $c>0$ for $c^np(z/c)$ see \cite{vdSl70},  Def.~1.7). Hence,  when estimating $\Gamma(p)/\rho(p)$ for non-trivial $p$ we may  assume (compare Theorem 2.6 in \cite{vdSl70}) that $\rho(p)=1$ and hence $\sum_{i=0}^{n-1} |a_i| =1.$ This implies that $ \sqrt[k]{|a_{n-k}|} \leq 1,$ and moreover $\tau \leq 1.$ Thus, the relative overestimation of the Cauchy bound by our new bound $\Gamma(p)$ does not exceed the factor $\sqrt{9.45} \sim 3.0741.$ 

 Finally, a root-modulus bound  simultaneously having good relative overestimation of the Cauchy bound $ \rho(p)$ and the maximum root-modulus $ \mu(p)$  can be defined by  $\min \{\Gamma(p); F(p)\}.$  This bound is at most double the Cauchy bound (like $F(p)$, cf. \cite{vdSl70}, Th. 2.6), and retains the very good relative overestimation of $\Gamma(p)$ (see Theorem \ref{QualityNewBound} above) at the same asymptotic complexity as Fujiwara's bound $F(p).$\\
 
 \subsection{Outlook} 
 While the above result brings a certain closure to the search for measurably good, \emph{absolute root bounds} for root-moduli  a lot remains to be investigated. A  new, different quality measure could be established via a benchmark suite of polynomials with several, wide ranging coefficient and zero distributions.   It would be valuable to determine bounds, composed via algebraic functions of coefficients, which improve over \eqref{OPTIMUM}.  A study of such \emph{non-absolute} bounds should relate the computational effort to the quality of the bound. 
  

\bibliographystyle{plain}
{ 
	 }
\bibliographystyle{plain}
\end{document}